\documentclass[reqno,10pt]{amsart}
\newtheorem{Pa}{Paper}[section]
\newtheorem{Tm}[Pa]{{\bf Theorem}}
\newtheorem{La}[Pa]{{\bf Lemma}}

\newtheorem{Rk}[Pa]{{\bf Remark}}

\newtheorem{Dn}[Pa]{{\bf Definition}}
\newcommand{\R}{\mathbb R}
\newcommand{\raro}{\rightarrow}

\newcommand{\doubcolb}[4]{\left(\begin{array}{cc}#1&#2\\#3&#4\end{array}\right)}

\newcommand{\Gt}{\Theta}

\newcommand{\singcolb}[2]{\left(\begin{array}{c}#1\\#2\end{array}\right)}
\newcommand{\ben}{\begin{enumerate}}
\newcommand{\een}{\end{enumerate}}

\newcommand{\C}{\mathbb C}

\begin{document}
\newcommand{\eqd}{\buildrel\rm def\over =}

\title[Degenerate Hamburger moment problem]
{On degenerate Hamburger moment problem and extensions
of positive semidefinite Hankel block matrices}
\author{Vladimir Bolotnikov}
\address{Department of Mathematics,
The College of William and Mary,
Williamsburg VA 23187-8795, USA}

\begin{abstract}
   In this paper we consider two related objects: singular positive 
semidefinite Hankel block--matrices and associated degenerate truncated 
matrix Hamburger moment problems. The description of all solutions of a degenerate 
matrix Hamburger moment problem is given in terms of a linear fractional 
transformation. The case of interest is the Hamburger moment problem
whose Hankel block--matrix admits a  positive semidefinite Hankel  extension.
\end{abstract}
\maketitle 

This is the corrected version of the original paper \cite{bol}. The work was 
inspired by V. Dubovoj's paper \cite{dub1} containing the first systematic study
of degenerate matricial interpolation problems. Another sourse of inspiration 
must have been the paper by R. Curto and L. Fialkow \cite{cur} but  I was not 
aware of it then. The original paper contained several erratae and the author 
is very grateful to  A. Ben-Artzi and H. Woerdeman for indicating them. 
A short proof in Section 5  fixes these incorrectnesses. The remaining four sections
are mostly the same as in \cite{bol}.

   \section{Introduction}
\setcounter{equation}{0}
   The objective of this article is to describe the solutions of a
   degenerate truncated matrix Hamburger moment problem ${\bf HMP}$. 
   We start with a set of Hermitian matrices $s_{0},\ldots,s_{2n} \in
   \C^{m \times m}$ and let $K_{n}$ denote the Hankel block matrix
\begin{equation}
   K_{n}=(s_{i+j})^{n}_{i,j=0}.
\label{1.1}
\end{equation}
   Let ${\mathcal Z}(K_{n})$ denote the set of all solutions of the 
associated
   truncated Hamburger moment problem, i.e., the set of nondecreasing
   right continuous ${m \times m}$ matrix-valued functions
   $\sigma(\lambda)$ such that
\begin{equation}
   \int^{\infty}_{-\infty}\lambda^{k}d\sigma(\lambda)=s_{k} \hspace{10mm}
   (k=0,\ldots,2n-1)
\label{1.2}
\end{equation}
   and
\begin{equation}
   \int^{\infty}_{-\infty}\lambda^{2n}d\sigma(\lambda) \leq
   s_{2n}.
\label{1.3}
\end{equation}
   As in the scalar case (see [1: \S2.1]) ${\mathcal Z}(K_{n})$ is 
nonempty if and only
   if $K_{n}$ is positive semidefinite and, moreover, by a theorem of H. Hamburger and R. 
Nevanlinna [1: \S 3.1], the formula 
\begin{equation}
   w(z)=\int^{\infty}_{-\infty} \frac{d\sigma(\lambda)}{\lambda -z}
\label{1.4}   
\end{equation}
establishes  a one-to-one correspondence between ${\mathcal Z}(K_{n})$ 
and the class ${\mathcal R}(K_{n})$ of $\C^{m \times m}$--valued functions 
$w(z)$ analytic and with 
   positive semidefinite imaginary part in the  upper half plane $\C_{+}$ such 
that  uniformly in the angle
   $\{ z=\rho e^{i\theta}:\varepsilon
   \leq \theta \leq \pi-\varepsilon\;,\;\varepsilon>0\}$,
\begin{equation}
   \lim_{z \raro \infty}\left \{z^{2n+1}w(z)+\sum^{2n}_{k=0}
   s_{k} z^{2n-k}\right\} \geq 0.
\label{1.35}
\end{equation}
This correspondence reduces the ${\bf HMP}$ problem to a boundary interpolation problem 
of finding all $\C^{m \times m}$--valued Pick functions
   $w$ (which by definition are analytic and with positive semidefinite 
imaginary part in $\C_{+}$) with prescribed asymptotic behavour (\ref{1.35}) at infinity. 

   In this paper we follow the Potapov's method of the fundamental
   matrix inequality \cite{kov}. The starting point is the following
   theorem which describes the set ${\mathcal R}(K_{n})$ in terms of a 
   matrix inequality (see \cite[\S 1]{kov} for the proof).
\begin{Tm}
   Let $w$ be a $\C^{m \times m}$--valued function analytic in
   $\C_{+}$. Then $w$ belongs to $R(K_{n})$ if and only if it
   satisfies the inequality
\begin{equation}
   \left( \begin{array}{cc}
   K_{n} & (I-zF_{m,n})^{-1}(Uw(z)+M) \\  \\
   (w(z)^*U^*+M^*)(I-\bar{z}F_{m,n}^*)^{-1} & {\displaystyle\frac{w(z)-w(z)^{*}}{z-\bar{z}}}  
\end{array} \right)    \geq 0
\label{1.5}
\end{equation}
for every $z\in \C_{+}$, where
\begin{equation}
    F_{m,n}=    \left( \begin{array}{ccccc}
   0_{m} & & \ldots & & 0 \\ I_{m} & \ddots & & & \\
   0 &I_{m} &  & &\vdots \\ \vdots& \ddots & \ddots & \ddots & \\
   0 & \ldots & 0 & I_{m} & 0_{m} \end{array} \right)\in
   \C^{m(n+1) \times m(n+1)}
\label{1.6}
\end{equation}
is the matrix of the m-dimensional shift in
   $\C^{m(n+1)}$ and where $U, M \in \C^{m(n+1) \times m}$ are 
given by
\begin{equation}
   U=  \left( \begin{array}{c}
   I_{m} \\ 0 \\ \vdots \\ 0 \end{array} \right),
   \hspace {10mm} M=F_{m,n}K_{n}U=
   \left( \begin{array}{c} 0 \\ s_{0} \\ \vdots \\ s_{n-1}
   \end{array} \right).
\label{1.7}
\end{equation}
\label{T:1.1}
\end{Tm}
The matrix $K_{n}$ (the so--called {\it Pick matrix} of the ${\bf HMP}$)
   satisfies the following Lyapunov identity
\begin{equation}
   F_{m,n}K_{n} - K_{n}F_{m,n}^*= MU^{*}-UM^{*}
\label{1.8}
\end{equation}
   which can be easily verified with help of (\ref{1.1}), (\ref{1.6}) and (\ref{1.7}).

\smallskip

   The ${\bf HMP}$ is called {\em nondegenerate} if its Pick matrix $K_{n}$ is 
strictly positive and it is termed {\em degenerate} if 
$K_{n}$ is singular and positive semidefinite.
   The parametrization of all solutions to the inequality (\ref{1.5}) for the case 
   $K_n>0$ was obtained in \cite{kov} and will be recalled in Theorem \ref{T:1.3} below.
To formulate this theorem we first introduce some needed definitions and notations. We will 
denote bt ${\bf W}$ the class of $\C^{2m \times 2m}$--valued 
meromorphic functions $\Theta$ which are $J$--unitary on $\R$ and $J$--expansive in 
$\C_{+}$: 
\begin{equation}
    \Theta (z) J\Theta (z)^* =J \  \ \
    (z \in \R),   \, \; \; \ \ \ \ \
    \Theta (z) J\Theta (z)^* \geq J \  \ \; (z \in \C_+)
\label{1.85}
\end{equation}
    where 
\begin{equation}
    J= \left( \begin{array}{cc}
    0 & iI_{m} \\ -iI_{m} & 0 \end{array} \right).
\label{1.9}
\end{equation}
\begin{Dn}
A pair $\{p, \; q\}$ of $\C^{m \times m}$-valued functions meromorphic in $\C\backslash 
\R$ is called a {\it Nevanlinna pair} if
\begin{equation}
    \begin{array}{ll}
    (i) \; & {\rm det} \; (p(z)^*p(z)+q(z)^*q(z)) \not\equiv 0 
    \hspace{13mm}\mbox{(the nondegeneracy of the pair)} \\ \\
    (ii) \; & {\displaystyle\frac{q(z)^*p(z)-p(z)^*q(z)}{z-\bar{z}}}=
    \left(p(z)^{*},\hspace{1mm}q(z)^{*}\right) {\displaystyle\frac{J}{i(\bar{z}-z)}}
    \left( \begin{array}{c} p(z) \\ q(z)\end{array} \right) \geq 0
    \hspace{10mm} (\Im z\neq 0).\end{array}
\label{1.10}
\end{equation}
\label{D:1.3}
\end{Dn}
    A pair $\{p, \; q\}$ is said to be {\it equivalent} to the pair  
    $\{ p_{1}, \; q_{1} \}$ if there exists a $\C^{m \times m}$-valued
    function $\Omega$  $(det \; \Omega (z) \not\equiv 0)$ meromorphic in
    $\C\backslash \R$ such that $\; \; p_{1} = p \; \Omega \;$
    and $\; q_{1} = q \Omega $. The set of all $m\times m$ matrix valued Nevanlinna pairs 
    will be denoted by ${\bf N}_{m}$. 
\begin{Tm}
Let $K_n$ be a strictly positive matrix given by $(\ref{1.1})$ and let $F_{m,n}, \; U$ and 
$M$ be defined by $(\ref{1.6})$, $(\ref{1.7})$. Then
\begin{enumerate}
\item The function 
\begin{equation}
    \Theta(z)=\doubcolb{\theta_{11}(z)}{\theta_{12}(z)}{\theta_{21}(z)}{\theta_{22}(z)}
    =I_{2m} + z \left( \begin{array}{c} M^{*} \\ -U^{*}
    \end{array} \right) (I-zF_{m,n}^{*})^{-1}K_{n}^{-1}(U,M)
\label{1.11}
\end{equation} 
belongs to the class ${\bf W}$.
\item The formula
\begin{equation}
    w(z)=(\theta_{11}(z)p(z)+\theta_{12}(z)q(z))
    (\theta_{21}(z)p(z)+\theta_{22}(z)q(z))^{-1}
\label{1.12}
\end{equation}
gives all the solutions $w$ to the inequality  $(\ref{1.5})$ when
 $\{p, \; q\}$ varies in ${\bf N}_{m}$.
\item Two pairs $\{p(z), \; q(z)\}$ and $\{p_{1}(z), \; q_{1}(z)\}$
    lead by $(\ref{1.12})$ to the same function $w(z)$ if and only if these pairs are
    equivalent.
\end{enumerate}
\label{T:1.3}
\end{Tm}
The degenerate scalar ${\bf HMP}$ is simple: 
${\mathcal R}(K_{n})$ 
    consists of the unique rational function $w(z)$ (this follows immediately 
    from (\ref{1.5})). In the degenerate
    {\it matrix} case, the description of ${\mathcal R}(K_n)$ depends on 
the degeneracy of $K_n$, but we still have a
    parametrization of all the solutions as a linear fractional
    transformation (\ref{1.12}) with the coefficient matrix $\Theta$ from the class ${\bf 
    W}$
    and for a suitable choice of parameters $\{p, \; q\}$ (see Theorem \ref{T:4.6} below).
    To construct the coefficient matrix of the degenerate ${\bf HMP}$, we follow the method 
    of V. Dubovoj which was applied in [4] to the degenerate Schur problem.
    Note that if $\det \theta_{22}\not \equiv 0$, the transformation (\ref{1.12}) can be
    written as
\begin{equation}
     w(z)=\psi_{11}(z)+\psi_{12}(z)p(z)(\psi_{22}(z)p(z)+q(z))^{-1}\psi_{21}
\label{1.13}
\end{equation}
    where
\begin{equation}
    \psi_{11}=\theta_{11}\theta_{22}^{-1}, \; \; \; \;
    \psi_{12}=\theta_{11}-\theta_{12}\theta_{22}^{-1}\theta_{22}, \; \; \; \;
    \psi_{21}=\theta_{22}^{-1}, \; \; \; \;
    \psi_{22}=\theta_{22}^{-1}\theta_{21}
\label{1.14}
\end{equation}
    and it turns out that the function 
    $\Psi(z)={\scriptsize\doubcolb{\psi_{11}(z)}{\psi_{12}(z)}{\psi_{21}(z)}{\psi_{22}(z)}}$
    is a Pick function (i.e. analytic and with positive semidefinite 
imaginary part in $\C_+$).
    If $\det \theta_{22}\equiv 0$, formulas (\ref{1.14}) make no sense, but nevertheless the
    set ${\mathcal R}(K_{n})$ can be parametrized by the transformation 
(\ref{1.13}) with a 
    coefficient matrix $\Psi$ from the Pick class. This $\Psi$ can be constructed
as a characteristic
    function of certain unitary colligation associated with the initial data $\{s_j\}$ of 
the problem. This approach (see \cite{kh1}) is much more stable with respect to a possible 
    degeneracy  of the Pick matrix $K_n$. The degenerate ${\bf HMP}$ will be discussed in 
some more detail in Section 2. 

    \section {Positive semidefinite Hankel extensions of Hankel block matrices}
   \setcounter{equation}{0}

Let ${\bf {\mathcal H}}_{m,n}$ be the set of all
positive semidefinite Hankel block matrices of the form (\ref{1.1}). We say that 
a matrix $K_{n}\in{\bf {\mathcal H}}_{m,n}$  {\it admits a positive semidefinite Hankel 
extension} if there exist Hermitian matrices $s_{2n+1},s_{2n+2} \in \C^{m \times m}$ such 
that the block matrix $K_{n+1}=(s_{i+j})^{n+1}_{i,j=0}$ is still positive 
semidefinite. The class of such matrices will be denoted by ${\bf {\mathcal H}}_{m,n}^{+}$:
\begin{equation}
{\bf {\mathcal H}}_{m,n}^{+}=\left\{K_n\in{\bf {\mathcal H}}_{m,n}: \, 
(s_{i+j})^{n+1}_{i,j=0}\ge 0 \; \mbox{for some}\; s_1=s_1^*\; \mbox{and}\; 
s_2=s_2^*\right\}.
\label{def1}
\end{equation}
In the scalar case $(m=1)$ every positive semidefinite Hankel matrix
admits a positive semidefinite Hankel extension
and therefore, ${\bf {\mathcal H}}_{1,n}^{+}={\bf {\mathcal H}}_{1,n}$. For 
$n\ge 2$,  ${\bf {\mathcal H}}_{m,n}^{+}$ is a proper subset of ${\bf {\mathcal H}}_{m,n}$
as can be seen from the example
$$
    K_2=\left( \begin{array}{cc} s_0& s_1 \\ s_1 & s_2 \end{array} \right),\qquad
    s_0=s_1=\left( \begin{array}{cc}
    1& 0 \\ 0 & 0 \end{array} \right)\quad{\rm and}\quad
    s_2=\left( \begin{array}{cc} 1& 0 \\ 0 & 1 \end{array} \right).
$$
We introduce two more subsets
of ${\bf {\mathcal H}}_{m,n}$:
\begin{equation}
\widetilde{\mathcal H}_{m,n}:=\left\{K_n\in{\bf {\mathcal H}}_{m,n}: \, {\bf P}_{Ker K_{n-1}} 
{\scriptsize \left( \begin{array}{c}
   s_{n+1} \\ \vdots \\ s_{2n} \end{array} \right)} =0\right\}
\label{dop}
\end{equation}
and 
\begin{equation}
\widehat{\mathcal H}_{m,n}:=\left\{K_n\in{\bf {\mathcal H}}_{m,n}: \, s_{2n}=\int_{-\infty}^{\infty}
\lambda^{2n}d\sigma(\lambda) \; \mbox{for some}\; \sigma\in{\mathcal Z}(K_n)\right\}.
\label{def2}
\end{equation}
Thus, $\widehat{\mathcal H}_{m,n}$ consists of all matrices $K_n\in{\bf {\mathcal H}}_{m,n}$, the 
associated truncated Hamburger moment problem admits an ``exact'' solution $\sigma$ such 
that 
\begin{equation}
   \int^{\infty}_{-\infty}\lambda^{k}d\sigma(\lambda)=s_{k} \hspace{10mm}
   (k=0,\ldots,2n),
\label{2.18}
\end{equation}
that is, with equality for the last assigned moment $s_{2n}$ rather than inequality 
\eqref{1.3}. In \eqref{dop} and in what 
follows, ${\bf P}_{Ker K}$ denotes  the orthogonal projection onto the kernel of $K$. We will show 
below that 
\begin{equation}
{\bf {\mathcal H}}_{m,n}^{+}=\widetilde{\mathcal H}_{m,n}=\widehat{\mathcal H}_{m,n}
\label{dop1}   
\end{equation}
which will provide therefore, several equivalent characterizations of Hankel block matrices 
admitting positive semidefinite Hankel extensions. The following two propositions can be easily 
verified.
\begin{La}
The block matrix
   $T=(t_{ij})^{n}_{i,j=0}\hspace{2mm}(t_{ij} \in \C^{r\times l})$
   is Hankel if and only if
\begin{equation}
   F^{*}_{l,n}(F_{l,n}T-TF^{*}_{r,n})F_{r,n}=0
\label{2.1}
\end{equation}
   where $F$ is a shift matrix defined via $(\ref{1.6})$.
\label{L:2.1}
\end{La}
\begin{La}
   Let $K,V \in \C^{N\times N}$ and $A \in \C^{N\times r}$ 
   be matrices such that $K=K^{*}$ and $\det \ V\neq 0$. Then,
   ${\bf P}_{Ker K}A={\bf P}_{Ker VKV^{*}}VA$.
\label{L:2.2}
\end{La}
Given a $K\ge 0$, let $Q$ be a matrix such that
\begin{equation}
   QKQ^{*}>0\hspace{5mm} {\rm and}\hspace{5mm}{\rm rank} \ QKQ^{*}={\rm rank} \ K.
\label{2.2}
\end{equation}
   We define the pseudoinverse matrix $K^{[-1]}$ by
\begin{equation}
   K^{[-1]} = Q^{*}\left(QKQ^{*}\right)^{-1}Q.
\label{2.3}
\end{equation}
 Since the pseudoinverse matrix depends on the choice of $Q$, it is not uniquely
 defined. 
\begin{La}
For every choice of $K^{[-1]}$,
\begin{equation}
   I-KK^{[-1]}=\left(I-KK^{[-1]}\right){\bf P}_{KerK}.
\label{2.4}
\end{equation}
\label{L:2.3}
\end{La}
{\bf Proof:} By (\ref{2.2}), every vector $f$
   can be decomposed as $\hspace{2mm}f=g+hQ\hspace{2mm}$ for some\\
   $g\in Ker \ K$ and $h\in \C^{1\times {\rm rank}K}$. Therefore,
   $$
   f\left(I-KK^{[-1]}\right)=
   (g+hQ)\left(I-KQ^{*}(QKQ^{*})^{-1}Q\right)=g
   $$
   which implies (\ref{2.4}).\qed 
\begin{La}
The block matrix
   ${\scriptsize \left( \begin{array}{cc}
   K & B \\ B^{*} & C \end{array} \right)}$ is positive semidefinite if 
and only if 
   $$
   K\geq 0,\hspace{4mm}{\bf P}_{kerK}B=0\hspace{3mm} {\rm and} \hspace{3mm}
   R= C-B^{*}K^{[-1]}B\geq 0.
   $$
   Moreover, if ${\scriptsize \left( \begin{array}{cc}
   K & B \\ B^{*} & C \end{array} \right)}\geq 0$, then the matrix $R$
   does not depend on the choice of $K^{[-1]}$.
\label{L:2.4}
\end{La}
   {\bf Proof:} The first assertion of lemma follows from the
   factorization
   $$
   \doubcolb{K}{B}{B^{*}}{C}=\doubcolb{I}{0}{B^{*}K^{[-1]}}{I}
   \doubcolb{K}{0}{0}{R}\doubcolb{I}{K^{[-1]}B}{0}{I}
   $$
   which in view of (\ref{2.4}), is valid if and only if ${\bf P}_{kerK}B=0$.\\
   Furthermore, let $C$ admit two different representations
   $C=R_i+B^{*}K_i^{[-1]}B\hspace{2mm}(i=1,2)$. Then
\begin{equation}
   R_1-R_2=B^{*}\left( K_2^{[-1]}-K_1^{[-1]}\right)B.
\label{2.5}
\end{equation}
   In view of (\ref{2.4}),
   $$
   K\left( K_2^{[-1]}-K_1^{[-1]}\right)B=\left\{\left(I-KK_1^{[-1]}\right)-
   \left(I-KK_2^{[-1]}\right)\right\}{\bf P}_{KerK}B=0.
   $$
   Since ${\scriptsize \left( \begin{array}{cc} K & B \\ B^{*} & C \end{array}
   \right)}\geq 0$, then also $B^{*}\left( K_2^{[-1]}-K_1^{[-1]}\right)B=0$
   which both with (\ref{2.5}) implies $R_1=R_2$. \qed 
\begin{La}
 Let  $K_{n} \in {\bf {\mathcal H}}_{m,n}$ and let ${\mathcal L}$ be the
subspace of $\C^{1 \times m}$ given by
\begin{equation}
   {\mathcal L}=\{f \in \C^{1 \times m}: (f_{0},\ldots,f_{n-2},f) \in
   KerK_{n-1}
   \hspace{2mm}for\hspace{2mm}some\hspace{2mm} f_{0},\ldots,f_{n-2}
   \in \C^{1\times m}\}.
\label{2.6}   
\end{equation}
Then $K_n$ belongs to $\widetilde{\mathcal H}_{m,n}$, that is 
\begin{equation}
{\bf P}_{Ker K_{n-1}} {\scriptsize \left( \begin{array}{c}
   s_{n+1} \\ \vdots \\ s_{2n} \end{array} \right)} =0,
\label{2.7}
\end{equation}
if and only if the block $s_{2n}$ is of the form
\begin{equation}
   s_{2n}=(s_{n},\ldots,s_{2n-1})K_{n-1}^{[-1]}(s_{n},\ldots,s_{2n-1})^{*}+R
\label{2.8}
\end{equation}
   for some positive semidefinite matrix $R\in \C^{m \times m}$
   which vanishes on the subspace ${\mathcal L}$ and does not depend on
the choice of $K_{n-1}^{[-1]}$.
\label{L:2.5}
\end{La}
{\bf Proof:} Since $K_n\geq 0$, then by Lemma 2.4,
$$
s_{2n}-(s_{n},\ldots,s_{2n-1})K_{n-1}^{[-1]}(s_{n},\ldots,s_{2n-1})^{*}\geq 0
$$
and therefore, $s_{2n}$ admits a representation (\ref{2.8}) 
for some $R\geq 0$. Moreover, since $K_{n}\geq 0$, then for
   every vector $(f_0,\ldots,f_{n-1})$ from $KerK_{n-1}$
   $$
   (f_0,\ldots,f_{n-1})
   \left( \begin{array}{ccc} s_{1} & \ldots &
   s_{n} \\ \vdots & &  \vdots \\ s_{n} &
   \ldots & s_{2n-1}\end{array} \right) =0
   $$
   and therefore,
\begin{eqnarray}
   f_{n-1}(s_{n},\ldots,s_{2n-1}) &=&-(f_0,\ldots,f_{n-2})
   \left( \begin{array}{ccc} s_{1} & \ldots &
   s_{n} \\ \vdots & & \vdots \\ s_{n-1} &
   \ldots & s_{2n-2}\end{array} \right)\nonumber \\
  & =&-(0,f_0,\ldots,f_{n-2}) K_{n-1}.\label{dop3}
\end{eqnarray}
Thus, 
\begin{eqnarray}
&& f_{0}s_{n+1} + \ldots +f_{n-2}s_{2n-1} + f_{n-1}(s_{n},\ldots,s_{2n-1})
   K_{n-1}^{[-1]}(s_{n},\ldots,s_{2n-1})^{*}\nonumber\\
&&=(0,f_0,\ldots,f_{n-2})\left \{ I-K_{n-1}K_{n-1}^{[-1]}\right \}
(s_{n},\ldots,s_{2n-1})^{*}\nonumber\\
&&=(0,f_0,\ldots,f_{n-2})\left \{ I-K_{n-1}K_{n-1}^{[-1]}\right \}
   {\bf P}_{Ker K_{n-1}}(s_{n},\ldots,s_{2n-1})^{*} =0\label{dop4}
\end{eqnarray}
where the first equality holds due to \eqref{dop3}, the second follows by
(\ref{2.4}) and the last one holds since $K_n\geq 0$ and therefore, ${\bf P}_{Ker 
K_{n-1}}(s_{n},\ldots,s_{2n-1})^{*} =0$. Comparing \eqref{dop4} with (\ref{2.8})
gives
\begin{equation}
   f_{0}s_{n+1} + \ldots +f_{n-1}s_{2n}=f_{n-1}R. 
\label{2.13}
\end{equation}
It remains to show that $R$ vanishes on the subspace ${\mathcal L}$ if and only if 
\eqref{2.7} holds. To this end, let us observe that condition (\ref{2.7}) means that
   $f_{0}s_{n+1} + \ldots +f_{n-1}s_{2n}=0$  for every vector
   $(f_{0},\ldots ,f_{n-1}) \in KerK_{n-1}.$
   The latter is equivalent, in view of (\ref{2.13}) and (\ref{2.6}), to
   $f_{n-1}R=0$ for all $f_{n-1}\in {\mathcal L}$. By Lemma 2.4, the
matrix
   $\hspace{3mm} R=
   s_{2n}-(s_{n},\ldots,s_{2n-1})K_{n-1}^{[-1]}(s_{n},\ldots,s_{2n-1})^{*}$
   does not depend on the choice of  $K_{n-1}^{[-1]}$.\qed 
\begin{La}
Let ${\mathcal H}_{m,n}^+$,  $\widetilde{\mathcal H}_{m,n}$ and $\widehat{\mathcal H}_{m,n}$
be the classes defined in \eqref{def1}--\eqref{def2}. Then
\begin{equation}
{\mathcal H}_{m,n}^+\subseteq \widehat{\mathcal H}_{m,n}\subseteq \widetilde{\mathcal H}_{m,n}.
\label{2.13a}
\end{equation}
\label{L:2.6}
\end{La}
{\bf Proof:}  Let $K_{n+1}$ be a positive semidefinite Hankel
   extension of $K_n$. Since $K_{n+1}\geq 0$, by the solvability criterion for the
associated Hamburger moment problem, the set ${\mathcal Z}(K_{n+1})$ is nonempty.
Furthermore, for every $\sigma\in {\mathcal Z}(K_{n+1})$ 
$$
   \int^{\infty}_{-\infty}\lambda^{k}d\sigma(\lambda)=s_{k}\quad
   (k=0,\ldots,2n+1)\mbox{and}\quad \int^{\infty}_{-\infty}\lambda^{2n+2}d\sigma(\lambda) \leq
   s_{2n+2}
$$
and therefore, $K_n\in\widehat{\mathcal H}_{m,n}$ which proves the first containment in 
\eqref{2.13a}.

Now let us assume that $K_n$ belongs to $\widehat{\mathcal H}_{m,n}$ and let 
$d\sigma$ be the measure satisfying conditions  (\ref{2.18}). Then 
\begin{equation}
   K_n=\int^{\infty}_{-\infty}\left(I_m,\ldots,\lambda ^{n}I_m\right)^{*}
   d\sigma(\lambda)\left(I_m,\ldots,\lambda ^{n}I_m\right). 
\label{2.19}
\end{equation}
   Let $\hspace{2mm}{\bf f}=(f_0,\ldots,f_{n-1})\in \C^{1\times mn}\hspace{2mm}$
   be a vector from $Ker\hspace{1mm}K_{n-1}$. Then
   $\hspace{2mm}\int^{\infty}_{-\infty}f(\lambda)d\sigma(\lambda)f(\lambda)^{*}=0,
   \hspace{2mm}$ where
\begin{equation}
   f(\lambda)=f_0+\lambda f_1+\ldots +\lambda ^{n-1}f_{n-1}={\bf f}
   \left(I_m,\ldots,\lambda ^{n-1}I_m\right)^{*}.
\label{2.20}
\end{equation}
   In particular, for every choice of $\hspace{2mm}-\infty<a<b<+\infty$,
\begin{equation}
   \int^{b}_{a}f(\lambda)d\sigma(\lambda)f(\lambda)^{*}=0.
\label{2.21}
\end{equation}
   Let $\; g\in \C^{1\times m}\hspace{1mm}$ be an arbitrary nonzero
   vector.  By the Cauchy inequality,
   $$
   \int^{b}_{a}f(\lambda)d\sigma(\lambda)\lambda ^{n+1}g^{*}\hspace{2mm}\leq
   \left( \int^{b}_{a}f(\lambda)d\sigma(\lambda)f(\lambda)^{*}\hspace{1mm}
   \int^{b}_{a}\lambda ^{2n+2}gd\sigma(\lambda)g^{*}\right)^{\frac{1}{2}}
   $$
   which in view of (\ref{2.21}) implies $\hspace{2mm}\int^{b}_{a}f(\lambda)
   d\sigma(\lambda)\lambda ^{n+1}g^{*}=0$. Since  $\; a,b \in \R \;$ and 
   $\; g\in \C^{1\times m} \;$ are arbitrary, then
   $$
   \int^{\infty}_{-\infty}f(\lambda)d\sigma(\lambda)\lambda ^{n+1}I_m=0
   $$
   which on account of (\ref{2.18})--(\ref{2.20}) can be rewritten as
\begin{equation}
   {\bf f}(s_{n+1},\ldots,s_{2n})^{*}=0.
\label{2.22}
\end{equation}
   Thus, every vector $\; {\bf f}\in Ker \ K_{n-1} \;$ satisfies (\ref{2.22}) 
   or in other words,   
   $\; {\bf P}_{Ker K_{n-1}} \left(s_{n+1},\ldots,s_{2n} \right)^{*} =0 \;$
   and therefore, $K_{n} \in\widetilde{\mathcal H}_{m,n}$, which completes the proof of the 
second inclusion in \eqref{2.13a}.\qed 

\medskip

   In connection with the last lemma we consider the following question:
   {\it to describe all matrices $s\in \C^{m\times m}$ such that
   $s=\int^{\infty}_{-\infty}\lambda^{2n}d\sigma(\lambda)$ for some
   $\sigma\in {\mathcal Z}(K_n)$.}    
\begin{La}
Let $K_n\ge 0$ be a block matrix of the form \eqref{1.1} with the block $s_{2n}$
of the form 
\begin{equation}
   s_{2n}=\left(s_{n},\ldots,s_{2n-1}\right)K_{n-1}^{[-1]}
   \left(s_{n},\ldots,s_{2n-1} \right)^{*} +R
\label{2.23}
\end{equation}
for some matrix $R\ge 0$ (which does not depend on the choice of $K_{n-1}^{[-1]}$)
and let $s\in \C^{m\times m}$ be defined by \begin{equation}
   s=\int^{\infty}_{-\infty}\lambda^{2n}d\sigma(\lambda)
\label{2.24}  
\end{equation}
for some $\sigma\in {\mathcal Z}(K_n)$. Then  there exists a positive semidefinite matrix 
$R_{0}\leq R$ which  vanishes on the subspace ${\mathcal L}$ defined by $(\ref{2.6})$ and 
such that
\begin{equation}
   s=(s_{n},\ldots,s_{2n-1})K_{n-1}^{[-1]}            
   \left(s_{n},\ldots,s_{2n-1} \right)^{*}+R_{0} \quad( 0\le R_0\le R\quad\mbox{and}\quad 
R_0\vert_{\mathcal L}=0).
\label{2.25}
\end{equation}
\label{L:2.7}
\end{La}
   {\bf Proof:} Let $s$ be of the form (\ref{2.24}) for some  $\sigma\in
   {\mathcal Z}(K_n)$. We introduce the Hankel block matrix
\begin{equation}
   \widetilde{K}_{n}= \left( \begin{array}{cccc} s_{0} & \ldots & s_{n-1} & s_{n} \\
   \vdots & & & \vdots \\ s_{n-1} & & &
   s_{2n-1} \\ s_{n} & \ldots & s_{2n-1} & s \end{array} \right)
\label{2.26}
\end{equation}
   which differs from $K_n$ only by the block  $\widetilde{s}_{2n}=s$.
Thus, $\widetilde{K}_{n} \in\widehat{\mathcal H}_{m,n}$. Therefore,  $\widetilde{K}_{n} 
\in\widetilde{\mathcal H}_{m,n}$, by Lemma \ref{L:2.6}. By Lemma \ref{L:2.5}, 
   the block  $\widetilde{s}_{2n}=s$ admits representation (\ref{2.25}) for some $R_{0}\ge 0$ 
   vanishing on  ${\mathcal L}$.  The inequality $R_{0} \leq R$ follows  from
   (\ref{1.3}) and (\ref{2.23})--(\ref{2.25}). \qed
\begin{La}
Let $K_n\in{\mathcal H}_{m,n}$ be of the form \eqref{1.1}, let ${\mathcal L}$ be the subspace
given by \eqref{2.6}, let $s_{2n}$, $s$ and $\widetilde{K}_n$ be matrices defined by 
$(\ref{2.23})$, 
   $(\ref{2.25})$  and $(\ref{2.26})$ respectively, and let the positive semidefinite 
   $R_{0}: \C^{m} \raro \C^{m}$ be defined by 
\begin{equation}
   R_{0}h=  \left \{  \begin{array}{cc} 0 & for \hspace{4mm} h \in
   {\mathcal L}, \\ Rh & for \hspace{4mm} h\in {\mathcal L}^{\perp }.
   \end{array}\right. 
\label{2.27}
\end{equation}
   Then the Hamburger moment problems associated with the sets of matrices
   $\{s_{0},\ldots,s_{2n-1},s_{2n}\}$ and $\{s_{0},\ldots,s_{2n-1},s\}$ have
   the same solutions: ${\bf {\mathcal Z}}(K_n) = {\bf {\mathcal 
Z}}(\widetilde{K}_{n})$.
\label{L:2.8}
\end{La}
   {\bf Proof:} Let $\sigma$ belong to  ${\mathcal Z}(K_{n})$.
   By Lemma \ref{L:2.7}, the matrix $\widehat{s}= \int^{\infty}_
   {-\infty}\lambda^{2n}d\sigma(\lambda)$ admits a representation (\ref{2.25}) with
   a positive semidefinite matrix $\widehat{R}_{0}\leq R$ vanishing on
   ${\mathcal L}$. In view of (\ref{2.27}), $\widehat{R}_{0}\leq R_{0}$.
   Therefore, $\widehat{s}\leq s$ and $\sigma\in  {\mathcal Z}(\widetilde{K}_n)$. 
   So,  ${\mathcal Z}(K_{n})\subseteq {\mathcal Z}(\widetilde{K}_n)$. The 
converse 
inclusion
   follows from the inequality $s\leq s_{2n}$. \qed
\begin{Rk} By Lemmas \ref{L:2.5} and \ref{L:2.8}, 
we can assume without loss of generality that the Pick 
matrix of the  ${\bf HMP}$ belongs to $\widetilde{\mathcal H}_{m,n}$. 
\label{R:2.9}
\end{Rk}
   Otherwise we replace the block $s_{2n}$ (which is necessarily of the
   form (\ref{2.23})) by the block $\widetilde{s}_{2n}=s$ defined by (\ref{2.25}), 
   (\ref{2.27}). By Lemma \ref{L:2.5}, 
   $\widetilde{K}_{n}\in \widetilde{\mathcal H}_{m,n}$ and we describe the set 
   ${\mathcal Z}(\widetilde{K}_{n})$ of solutions of this new moment problem, 
which  coincides, by Lemma \ref{L:2.8}, with ${\mathcal Z}(K_n)$.

   \section{ The coefficient matrix of the problem}
\setcounter{equation}{0}

   The coefficient matrix $\Gt$ of the nondegenerate ${\bf HMP}$ given by the 
   formula (\ref{1.11}) is the matrix polynomial of $deg \ \Gt = n+1$ and 
   (\ref{1.11}) is a realization of $\Theta$ with state space equal $\C^{m(n+1)}$.
   In this section we obtain some special decomposition (see formula (\ref{3.13}) below) of 
the state space which will allow us to construct the analogue of (\ref{1.11}) for $K_{n}$ not 
strictly positive (formula  (\ref{3.25})). The idea is simple: to replace in
   (\ref{1.11}) the inverse  of the matrix $K_n$ (which does not exist for the
   degenerate case) by its pseudoinverse. However after this replacement the  
   function $\Theta$ may lose its $J$--properties (\ref{1.85}) which are essential
   for the description (\ref{1.12}) to be in force. This suggests the following   
   question: is there exist a pseudoinverse  matrix $K_{n}^{[-1]}$ of the form 
   (\ref{2.3}) such that the function    
$$
   \Theta(z)=I_{2m} + z \left( \begin{array}{c} M^{*} \\ -U^{*}  
    \end{array} \right) (I-zF_{m,n}^{*})^{-1}K_{n}^{[-1]}(U, \; M)
   $$
   still belongs to the class ${\bf W}$? We show in Lemmas \ref{L:3.2} and \ref{L:3.3} below
that such a pseudoinverse exists if (and in fact, only if) the Pick matrix $K_n$ belongs to 
the class $\widetilde{\mathcal H}_{m,n}$. Recall that for the degenerate matricial Schur 
problem such a pseudoinverse  always exists (see \cite{dub1}).
\begin{La}
Let $T_{n}=(t_{i+j})^{n}_{i,j=0} \in
   \widetilde{\mathcal H}_{l,n} \; (t_i\in \C^{l\times l})$, let 
$t_0>0$ and let $\widehat{T}_{n-1}$ be the block matrix defined as
\begin{equation}
   \widehat{T}_{n-1}=D^{-1}_{n} \left \{ {\bf S}-{\mathcal T}_nt_{0}^{-1}
   {\mathcal T}_n^{*}\right\}D^{-*}_{n}
\label{3.1}
\end{equation}
   where
\begin{equation}  
   D_{n}= 
\left( \begin{array}{cccc} t_{0} & 0 &
   \ldots & 0 \\
   t_{1} & & &  \\ \vdots & \ddots & \ddots & 0
   \\ t_{n-1} & \ldots & t_{1} & t_{0}
   \end{array} \right),\quad{\bf S}=(t_{i+j})^{n}_{i,j=1},\hspace{6mm}
   {\mathcal T}_n=\left( \begin{array}{c}
   t_{1} \\ \vdots \\ t_{n} \end{array} \right).
\label{3.2}
\end{equation}
   Then $\widehat{T}_{n-1}$ is a Hankel block matrix: 
\begin{equation}  
   \widehat{T}_{n-1}=(\widehat{t}_{i+j})^{n-1}_{i,j=0}
\label{3.3}
\end{equation}
   and moreover, $\widehat{T}_{n-1} \in \widetilde{\mathcal H}_{l,n-1}$.
\label{L:3.1}
\end{La}
   {\bf Proof:} Let $F_{l,n-1}$
   be the matrix defined via formula (\ref{1.6})    and let
\begin{equation}  
   \widetilde{U}:=(I_l,0,\ldots,0)^{*}\in\C^{ln\times l}.
\label{3.4}
\end{equation}
   We begin with the identities
\begin{equation}  
   D_{n}F_{l,n-1}=F_{l,n-1}D_{n},\hspace{6mm}
   \widetilde{U}^{*}F_{l,n-1}=0,\hspace{6mm}
   D_{n}\widetilde{U}-F_{l,n-1}{\mathcal T}_n=\widetilde{U}t_{0}
\label{3.5}
\end{equation}
   and
   $$
   F_{l,n-1}({\bf S}-{\mathcal T}_nt_{0}^{-1}{\mathcal T}_n^{*})-
   ({\bf S}-{\mathcal T}_nt_{0}^{-1}{\mathcal T}_n^{*})F_{l,n-1}^{*}
   ={\mathcal T}_nt_{0}^{-1}\widetilde{U}^{*}D_{n}^{*}
   -D_{n}\widetilde{U}t_{0}^{-1}{\mathcal T}_n^{*}
   $$
   which follow immediately from (\ref{1.6}), (\ref{3.2}) and (\ref{3.4}). Using 
these   identities we get
   $$
   \begin{array}{l}F_{l,n-1}^{*}(F_{l,n-1}\widehat{T}_{n-1}-
   \widehat{T}_{n-1}F_{l,n-1}^{*})F_{l,n-1} \\ \\
   =F_{l,n-1}^{*}D^{-1}_{n}\left\{F_{l,n-1}\left( {\bf S}-
   {\mathcal T}_nt_{0}^{-1} {\mathcal T}_n^{*} \right) -
   \left({\bf S}-{\mathcal T}_nt_{0}^{-1} 
   {\mathcal T}_n^{*}\right)F_{l,n-1}^{*}\right\}D^{-*}_{n}F_{l,n-1} \\ \\
   =F_{l,n-1}^{*}D^{-1}_{n}\left\{{\mathcal 
T}_nt_{0}^{-1}\widetilde{U}^{*}D_{n}^{*}
   -D_{n}\widetilde{U}t_{0}^{-1}{\mathcal T}_n^{*}\right\}D^{-*}_{n}F_{l,n-1} 
\\ \\
   =F_{l,n-1}^{*}D^{-1}_{n}{\mathcal T}_nt_{0}^{-1}\widetilde{U}^{*}F_{l,n-1}-
   F_{l,n-1}^{*}\widetilde{U}t_{0}^{-1}{\mathcal T}_n^{*}F_{l,n-1}=0 
\end{array}
   $$
and (\ref{3.3}) follows by Lemma \ref{L:2.1}. Since $D_{n}$ is invertible, the 
factorization formula
\begin{equation}
   T_{n}= \left( \begin{array}{cc} I_{l} & 0 \\
   {\mathcal T}_nt^{-1}_{0} & D_{n} \end{array} \right)
   \left( \begin{array}{cc}
   t_{0} & 0 \\ 0 & \widehat{T}_{n-1} \end{array} \right)
   \left( \begin{array}{cc} I_{l} & t^{-1}_{0}{\mathcal T}_n^{*} \\ 0 & 
D^{*}_{n}
   \end{array} \right)
\label{3.6}
\end{equation}
implies that $\widehat{T}_{n-1}\ge 0$ and thus, $\widehat{T}_{n-1} \in{\mathcal H}_{l,n-1}$.
It remains to verify that
\begin{equation}
   {\bf P}_{Ker \widehat {T}_{n-2}}{\scriptsize \left(\begin{array}{c}\widehat {t}_{n} \\
   \vdots \\ \widehat {t}_{2n-2}\end{array} \right)}=0.
\label{3.7}
\end{equation}
To this end, we first observe that
\begin{equation}
   {\bf P}_{KerT_{n-1}} \left({\mathcal T}_n, \; {\bf S}\right)=0
\label{3.8}
\end{equation}
since $T_{n}\ge 0$. Using the factorization of $T_{n-1}$ similar to (\ref{3.6}) we obtain
\begin{equation}
   \left( \begin{array}{cc} s_{0} & 0 \\ 0 & \widehat{T}_{n-2}   
   \end{array}\right)= \left( \begin{array}{cc} I & 0 \\
   -D^{-1}_{n-1}{\mathcal T}_{n-1}t^{-1}_{0} &
   D^{-1}_{n-1} \end{array} \right) T_{n-1} \left(
   \begin{array}{cc} I & -t^{-1}_{0}{\mathcal T}_{n-1}^{*}D^{-*}_{n-1} \\
   0 & D^{-*}_{n-1} \end{array} \right) 
\label{3.9}
\end{equation}
   where $\; D_{n-1}$ and ${\mathcal T}_{n-1} \;$ are defined via (\ref{3.2}).
Upon applying Lemma \ref{L:2.2} to the matrices 
$$
K=T_{n-1},\quad 
   V=\left( \begin{array}{cc} I & 0 \\ 
   -D^{-1}_{n-1}{\mathcal T}_{n-1} t^{-1}_{0} & D^{-1}_{n-1} \end{array} 
   \right)\quad\mbox{and}\quad A=\left({\mathcal T}_n, \; {\bf S}\right),
$$
and making use of (\ref{3.8}), (\ref{3.9}) we obtain
\begin{equation}
   {\bf P}_{Ker\widehat{T}_{n-2}}D^{-1}_{n-1}\left(- {\mathcal 
T}_{n-1}t^{-1}_{0},
   \; I_{mn}\right)\left({\mathcal T}_n, \; {\bf S}\right)=0.
\label{3.10}
\end{equation}
   From the block decomposition
   $D_{n}={\scriptsize \left( \begin{array}{cc} t_{0} & 0 \\
   {\mathcal T}_{n-1} & D_{n-1}\end{array} \right)}$ we have
\begin{equation}
   D^{-1}_{n}=\left( \begin{array}{cc} t^{-1}_{0} & 0 \\
   -D^{-1}_{n-1}{\mathcal T}_{n-1}t^{-1}_{0} & D^{-1}_{n-1} \end{array}
   \right).
\label{3.11}
\end{equation}
   Substituting (\ref{3.11}) into (\ref{3.1}) we obtain
   $$
   \left(
   \begin{array}{ccc} \widehat{t}_{1} & \ldots & \widehat{t}_{n}
   \\ \vdots & & \vdots \\ \widehat{t}_{n-1}
    & \ldots & \widehat{t}_{2n-2}\end{array} \right)
   = (0, \;  I_{m(n-1)})\widehat{T}_{n-1} 
   =D^{-1}_{n-1} \left(-{\mathcal T}_{n-1}t^{-1}_{0}, I_{mn}\right)\left( 
   {\bf S}-{\mathcal T}_nt^{-1}_{0}{\mathcal T}_n^{*} \right)D^{-*}_{n}.
   $$
   The last equality both with (\ref{3.10}) implies
   $$
   {\bf P}_{Ker \widehat {T}_{n-2}} \left( \begin{array}{ccc}
   \widehat{t}_{1} & \ldots
   & \widehat {t}_{n} \\ \vdots & & \vdots \\ \widehat{t}_{n-1} & \ldots &
   \widehat{t}_{2n-2}\end{array} \right)=0
   $$
   and, in particular, (\ref{3.7}), which completes the proof of lemma. \qed
\begin{La}
   Let   $K_n \in \widetilde{\mathcal H}_{m,n}$ and let 
   ${\rm rank}\hspace{1mm}K_{n}=r$. Then there exists $Q \in
   \C^{r\times (n+1)m}$ such that
\begin{equation}
   QK_{n}Q^{*}>0, \hspace{10mm} {\rm rank} \; QK_{n}Q^{*}=
   {\rm rank} \; K_{n}, \hspace{10mm}QF_{m,n}=NQ
\label{3.12}
\end{equation}
   for the shift $F_{m,n}$ defined by $(\ref{1.6})$ and some matrix 
$N \in
   \C^{r\times r}$. In other words, there exists a subspace
   $\; {\mathcal Q}={\rm Ran} \; Q \eqd \{y\in \C^{m(n+1)}: y= fQ$ \mbox 
{for some} $f \in \C^{r}\} \;$
   coinvariant with respect to $F_{m,n}$ and such that
\begin{equation}
   \C^{m(n+1)}=Ker \; K \; \dot{+} \;  {\mathcal Q}.
\label{3.13}
\end{equation}
\label{L:3.2}
\end{La}
   {\bf Proof:} We prove this lemma by induction. Let 
   $n=0$ and let $rank \ s_{0}=l\leq m$. Then there exists a
   unitary matrix ${\bf v}\in \C^{m\times m}$ such that 
\begin{equation}
   {\bf v}s_{0}{\bf v}^{*}= \left( \begin{array}{cl} t_{0} & 0 \\ 0 &
   0_{m-l}\end{array} \right) \hspace{10mm} (t_{0}>0),
\label{3.14}
\end{equation}
   and the matrix
\begin{equation}
   g=(I_{l},\hspace{2mm} 0){\bf v}\in \C^{l\times m}
\label{3.15}   
\end{equation}   
(considered as $Q$) clearly satisfies (\ref{3.12}).

\smallskip

   Let us suppose that the statement of the lemma holds for all integers up to $n-1$.
   Let as above, $\; rank \ s_{0}=l$ and let ${\bf v}$ and $g$  be matrices 
   defined by (\ref{3.14}), (\ref{3.15}). Since $K_{n}\in \widetilde{\mathcal
   H}_{m,n}$, we have $\; Ker \ s_{0}\subseteq Ker \ s_{i}$ for $i=1,\ldots,2n$, 
and then we have from (\ref{3.14}),
\begin{equation}
   {\bf v}s_{i}{\bf v}^{*}= \left( \begin{array}{cl} t_{i} & 0 \\ 0 &
   0_{m-l}\end{array} \right)\quad(t_{i}\in \C^{l\times l};
   \hspace{2mm}i=1,\ldots,2n).
\label{3.16}  
\end{equation}
In more detail, representations \eqref{3.16} for $i=1,\ldots,2n-1$ follow from 
positivity of $K_n$ along with its Hankel structure. Since $K_{n}$ belongs to  
$\widetilde{\mathcal H}_{m,n}$, equality \eqref{2.7} holds. Upon substituting 
decompositions \eqref{3.14} and \eqref{3.16} (for $i=1,\ldots,2n-1$) into 
\eqref{2.7}, one can easily see that $s_{2n}$ is necessarily of the form 
${\bf v}s_{2n}{\bf v}^{*}= \left( \begin{array}{cc} t_{2n} & \gamma \\ 0 &
 0\end{array} \right)$ for some $\gamma\in\C^{l\times(m-l)}$. Since 
$s_{2n}$ is Hermitian, $\gamma=0$ and representation \eqref{3.16} for $s_{2n}$ follows.

\smallskip

   From (\ref{3.14})--(\ref{3.16}) we obtain that 
   $\hspace{2mm}gs_{i}g^{*}=t_{i}\hspace{2mm}(i=0,\ldots,2n)\hspace{2mm}$ and
\begin{equation}
   T_n\equiv \left(t_{i+j}\right)^{n}_{i,j=0}=G_nK_nG_n^{*},\hspace{10mm}
   {\rm rank} \ T_n ={\rm rank} \ K_{n},
\label{3.17}   
\end{equation}   
   where $G_n$ is the $(n+1)l \times (n+1)m$ matrix defined by
\begin{equation}
   G_n = {\scriptsize \left( \begin{array}{ccc} g & & 0 \\  & \ddots &  \\ 0 &
   & g\end{array} \right)}.
\label{3.18}   
\end{equation}   
   Since $K_{n}\in \widetilde{\mathcal H}_{m,n}$, then it is readily seen that
   $T_{n}\in \widetilde{\mathcal H}_{l,n}$. Let 
$\hspace{2mm}\widehat{T}_{n-1},\hspace{2mm}D_{n}$ and 
${\mathcal T}_n$
   be matrices defined by (\ref{3.1}), (\ref{3.2}).
   Multiplying $K_{n}$ on the left by the matrix
\begin{equation}
   \Phi = \left( \begin{array}{cc} I_{l} & 0 \\ -D_{n}^{-1}{\mathcal T}_n
   t^{-1}_{0} & D_{n}^{-1}\end{array} \right)G_n 
\label{3.19}   
\end{equation}   
   and by $\Phi^{*}$ on the right we obtain, on account of of (\ref{3.18}) and (\ref{3.6}),
\begin{equation}
   \Phi K_{n}\Phi^{*}=\left( \begin{array}{cc} t_{0} & 0 \\ 0 &
   \widehat{T}_{n-1}\end{array} \right)  .
\label{3.20}   
\end{equation}   
   By Lemma \ref{L:3.1}, $\widehat{T}_{n-1} \in \widetilde{\mathcal H}_{l,n-1}$,
   and it follows from (\ref{3.17}), (\ref{3.19}) and (\ref{3.20})
   that $\; \; rank \ \widehat{T}_{n-1}= \ rank \ K_{n}-rank \ t_{0}=r-l. \;$
   Therefore, by the induction hypothesis, there exist matrices
   $\widetilde{Q}\in \C^{(r-l)\times ln}$ and $\widetilde{N}\in \C^{(r-l)\times
   (r-l)}$ such that
\begin{equation}
   \widetilde{Q}\widehat{T}_{n-1}\widetilde{Q}^{*}>0 \quad \mbox {and}  
   \quad\widetilde{Q}F_{l,n}=\widehat{N}\widehat{Q}.
\label{3.21}   
\end{equation}   
   We show that the matrices
\begin{equation}
   Q= \left( \begin{array}{cc} I_{l} & 0 \\ 0 & \widetilde{Q}\end{array}
   \right)\Phi\in \C^{r\times (n+1)},\hspace{12mm}
   N= \left( \begin{array}{cc}
   0_{l} & 0 \\ \widetilde{Q}\widetilde{U}t^{-1}_{0} 
    & \widetilde{N}\end{array} \right)\in \C^{r\times r}
\label{3.22}   
\end{equation}   
   (where $\; \widetilde{U} \;$ is the matrix given by (\ref{3.4}))
   satisfy (\ref{3.12}). Indeed, by (\ref{3.20})--(\ref{3.22},)
   $$
   QK_{n}Q^{*}= \left( \begin{array}{cc} I_{l} & 0 \\ 0 &
   \widetilde{Q}\end{array} \right)\left( \begin{array}{cc} t_{0} & 0
   \\ 0 & \widehat{T}_{n-1}\end{array} \right)\left( \begin{array}{cc}
   I_{l} & 0 \\ 0 & \widetilde{Q}^{*}\end{array} \right) >0
   $$
and
   $$
   {\rm rank} \ QK_{n}Q^{*}={\rm rank} \ t_{0}+
   {\rm rank} \ \widetilde{Q}\widehat{T}_{n-1}\widetilde{Q}^{*}=l+(r-l)={\rm rank} \  K_{n}.
   $$
We next make use of (\ref{3.19})--(\ref{3.21}) and of the block decompositions 
$$   
G_n=\doubcolb{g}{0}{0}{G_{n-1}}\quad\mbox{and}\quad
    F_{m,n}=\doubcolb{0}{0}{\widetilde{U}}{F_{m,n-1}} 
$$
to compute
$$
   QF_{m,n}= \doubcolb {0}{0}{\widetilde{Q}D_n^{-1}G_{n-1}\widetilde{U}}
   {\widetilde{Q}D_n^{-1}G_{n-1}F_{m,n-1}}=  
   \doubcolb {0}{0}{\widetilde{Q}D_n^{-1}\widetilde{U}g}{\widetilde{Q}D_n^{-1}
   F_{l,n-1}G_{n-1}}
$$
and 
\begin{eqnarray*}
   NQ & = &\doubcolb {0}{0}{(\widetilde{Q}\widetilde{U}-
   \widetilde{N}\widetilde{Q}D_n^{-1}{\mathcal T}_n)t_{0}^{-1}g}
   {\widetilde{N}\widetilde{Q}D_n^{-1}G_{n-1}} \\
   & =& \doubcolb {0}{0}{\widetilde{Q}D_n^{-1}(D_n\widetilde{U}-
   F_{l,n-1}{\mathcal T}_n) 
t_{0}^{-1}g}{\widetilde{Q}F_{l,n-1}D_n^{-1}G_{n-1}}.
\end{eqnarray*}
We now invoke equalities (\ref{3.5}) to verify that the 
right hand side matrices in the two last formulas coincide. Thus, 
$QF_{m,n}=NQ=0$, and the matrices  $Q$ and  $N$
   defined by (\ref{3.22}) satisfy (\ref{3.12}). This completes the proof. \qed 

\medskip

In what follows, the indeces will be omitted and 
by $K$ and $F$ we mean matrices $K_{n}$ and $F_{m,n}$ given
   by (\ref{1.1}) and (\ref{1.6}) respectively.
\begin{La}
Let $K\in\widetilde{\mathcal  H}_{m,n}$,
   let $Q$ be any matrix satisfying $(\ref{3.12})$ and let $F$, $U$, $M$,
   $J$ and  $K^{[-1]}$ be matrices given by $(\ref{1.6})$, $(\ref{1.7})$, 
   $(\ref{1.9})$ and $(\ref{2.3})$. Then the $\C^{2m\times  2m}$--valued function 
\begin{equation}
   \Theta(z)=I_{2m} + z \left( \begin{array}{c} M^{*} \\ -U^{*}
    \end{array} \right) (I-zF^{*})^{-1}K^{[-1]}(U, \; M)
\label{3.25}  
\end{equation}
   is of the class ${\bf W}$ and moreover,
\begin{equation}   
   \Theta(z)^{*}J\Theta(z)-J=i(\bar{z}-z)\left( \begin{array}{c}
    U^{*} \\ M^{*}\end{array} \right)K^{[-1]} (I-\bar{z}F)^{-1}K
   (I- zF^{*})^{-1}K^{[-1]}(U,M),
\label{3.26}  
\end{equation}
\begin{equation}   
   J-\Theta(z)^{-*}J\Theta^{-1}(z)=i(\bar{z}-z)\left(\begin{array}{c}
   U^{*} \\ M^{*}\end{array} \right)(I-\bar{z}F^{*})^{-1}
   K^{[-1]}(I- zF)^{-1}(U,M). 
\label{3.27}  
\end{equation}
\label{L:3.3}
\end{La}
Observe that the two first relations in (\ref{3.12}) enable us to construct
the pseudoinverse matrix $K^{[-1]}$ according to (\ref{2.3}) and the
third equality guarantees (\ref{3.26}) and (\ref{3.27}) to be in force.

\medskip

{\bf Proof:} Using (\ref{3.25}), \eqref{1.9} and (\ref{1.8}) we have
\begin{equation}
   \Theta(z)^{*}J\Theta(z)-J=i\left( \begin{array}{c}
    U^{*} \\ M^{*}\end{array} \right)L(z)(U, \; M)
\label{3.28}
\end{equation}
   where
\begin{equation}
   \begin{array}{ll}
   L(z) & = |z|^{2}K^{[-1]}(I-\bar{z}F)^{-1}\left\{MU^{*}-UM^{*} \right\}
   (I- zF^{*})^{-1}K^{[-1]} \\ \\
   &\hspace{4mm} +\bar{z}K^{[-1]}(I-\bar{z}F)^{-1}
   -z(I- zF^{*})^{-1}K^{[-1]}\\ \\
   & = (\bar{z}-z)K^{[-1]} (I-\bar{z}F)^{-1}K
   (I- zF^{*})^{-1}K^{[-1]} \\ \\
   &\hspace{4mm} +\bar{z}K^{[-1]}(I-\bar{z}F)^{-1}
   (I-KK^{[-1]})-z(I-K^{[-1]}K)(I- zF^{*})^{-1}K^{[-1]}.
   \end{array}
\label{3.29}  
\end{equation}
   It follows from (\ref{3.12}) that $\hspace{2mm}QF^j=N^jQ\hspace{3mm}(j=0,1,\ldots)$
   which both with (\ref{2.3}) implies
\begin{equation}
   K^{[-1]}F^j\left(I-KK^{[-1]}\right)=Q^{*}(QKQ^{*})^{-1}N^jQ
   \left(I-KQ(QKQ^{*})^{-1}Q \right)=0
\label{3.30}  
\end{equation}
   for $\hspace{2mm}j=0,1,\ldots$ Since
   $\hspace{2mm}(I-zF^{*})^{-1}=\sum_{j=0}^{n} z^{j}F^{*j}\hspace{2mm}$, then also
\begin{equation}    
   K^{[-1]}(I-zF)^{-1}\left(I-KK^{[-1]}\right)\hspace{20mm}(z\in\C).
\label{3.31}  
\end{equation}
   Substituting (\ref{3.31}) into (\ref{3.29}) and (\ref{3.29}) into (\ref{3.28}), we 
obtain (\ref{3.26}).
   Similarly,
\begin{equation}
   \Theta(z)J\Theta(z)^{*}-J=i\left( \begin{array}{c}
    M^{*} \\ -U^{*}\end{array} \right)(I- zF^{*})^{-1}\widetilde{L}(z)
    (I-\bar{z}F)^{-1}(M, \; -U)
\label{3.32}  
\end{equation}
   where
\begin{equation}
   \begin{array}{ll}
   \widetilde{L}(z) & =\bar{z}(I- zF^{*})K^{[-1]}-zK^{[-1]}(I-\bar{z}F)-
   |z|^{2}K^{[-1]}\left\{MU^{*}-UM^{*} \right\}K^{[-1]} \\ \\
   & = (\bar{z}-z)K^{[-1]} + |z|^{2}  K^{[-1]}F\left(I-KK^{[-1]}\right)
   -|z|^{2} \left(I-K^{[-1]}K\right)F^{*}K^{[-1]}. \end{array}
\label{3.33}  
\end{equation}
   Using (\ref{3.30}) for $\hspace{2mm}j=1 \;$ we obtain from (\ref{3.33}) that
   $\hspace{2mm}\widetilde{L}(z)=(\bar{z}-z)K^{[-1]} \;$ and by (\ref{3.32}),
\begin{equation}
   \Theta(z)J\Theta(z)^{*}-J=i(\bar{z}-z)\left( \begin{array}{c}
   M^{*} \\ -U^{*}\end{array} \right)(I- zF^{*})^{-1}
   K^{[-1]}(I-\bar{z}F)^{-1}(M, \; -U).
\label{3.34}  
\end{equation}
Relations \eqref{1.85} follow from \eqref{3.34} and thus, 
$\Theta\in {\bf W}$. Since it $\Theta$ is $J$--unitary on
   $\R$, then by the symmetry principle,\\   
   $\; \Theta^{-1} (z) = J \Theta(\bar{z})^* J \;$ which both  with 
   (\ref{3.34}) leads to
\begin{equation}
  \begin{array}{l}
   J - \Theta(z)^{-*} J\Theta^{-1}(z)
   =J( J-\Theta(\bar{z})J\Theta(\bar{z})^*)J    \\
   = i(z-\bar{z}) J\singcolb{M^*}{-U^{*}}(I-\bar{z}F^{*})^{-1}K^{[-1]}
   (I-zF)^{-1}(M,\hspace{1mm}-U)J   \end{array} 
\label{3.35}
\end{equation}
   and implies (\ref{3.27}).\qed
     
   \section{Parametrization of all solutions}
\setcounter{equation}{0}

   In this section we parametrize the set ${\mathcal R}(K_{n})$
   of all solutions of the degenerate ${\bf HMP}$ in terms of a
   linear fractional transformation. The following theorem can be found in \cite{kat, kov}.
\begin{Tm}
Let   $\Theta=\doubcolb {\theta_{11}}{\theta_{12}}{\theta_{21}}{\theta_{22}}$
   be the block decomposition of a $\C^{2m\times 2m}$--valued function $\Theta\in
   {\bf W}$ into four $\C^{m\times m}$--valued blocks. Then 
all $\C^{m\times m}$--valued analytic in $\C\backslash\R$ solutions $w$ to the inequality
\begin{equation}
   ( w(z)^*, \; I_m) \frac{\Theta(z)^{-*}J\Theta^{-1}(z)}{i(\bar{z}-z)}
   \singcolb{w(z)}{I_m} \geq 0
\label{4.1}
\end{equation}
   are parametrized by the formula 
\begin{equation}
   w(z)=(\theta_{11}(z)p(z)+\theta_{12}(z)q(z))
   (\theta_{21}(z)p(z)+\theta_{22}(z)q(z))^{-1}
\label{4.2}
\end{equation}
   when the parameter $\{p, \;  q\}$ varies in the set ${\bf N}_{m}$ of all 
   Nevanlinna pairs and satisfies
\begin{equation}
   {\rm det} \; (\theta_{21}(z)p(z)+\theta_{22}(z)q(z))\not\equiv 0;
\label{4.3}
\end{equation}
Moreover, two Nevanlinna pairs lead via $(\ref{4.2})$ to the same function $w$ if and only if 
   these pairs are equivalent.
\label{T:4.1}
\end{Tm}
\begin{La}
   Let $\{p, \; q\}\in {\bf N}_{m}$ be a Nevanlimma pair. Then  
\begin{equation}
   \det \; (p(z)+iq(z))\not\equiv 0, 
\label{4.4}
\end{equation}
   the function
\begin{equation}
   S(z)=(p(z)-iq(z))(p(z)+iq(z))^{-1}
\label{4.5}
\end{equation}
   is a $\C^{m\times m}$--valued contraction in $\C_{+}$ and moreover,
   two different pairs lead by $(\ref{4.5})$ to the same $s$ if and only if
   they are equivalent.
\label{L:4.2}
\end{La}
The proof is given in \cite{kat}. Observe that by (\ref{4.4}), every 
Nevanlinna pair $\{p, \; q\}$    satisfies the dual nondegeneracy property (compare with 
Definition \ref{D:1.3})
\begin{equation}
   {\rm det} \; (p(z)p(z)^*+q(z)q(z)^*) \not\equiv 0.
\label{4.6}
\end{equation}
\begin{La}
Let $\{p, \; q\}\in {\bf N}_{m}$ be a Nevanlimma
   pair   such that $\; \left(I_{\nu}, \; 0\right) \ p(z)\equiv 0 \; \; (\nu\leq m)$.
   Then $\{p, \; q\}$ is equivalent to a pair
\begin{equation}
   \left \{\left( \begin{array}{cc} 0_{\nu} & 0 \\ 0 &
   \widetilde{p}(z)\end{array}\right), \; \; \; \left( \begin{array}{cc} I_{\nu} & 0 \\ 0 &
   \widetilde{q}(z)\end{array}\right) \right\}\quad\mbox{for some}\quad \{\widetilde{p}, \; 
\widetilde{q}\}\in {\bf N}_{m-\nu}.
\label{4.65}
\end{equation}
\label{L:4.3}
\end{La}
   {\bf Proof:} By the assumption assumption, $p$ and $q$ are of the form 
\begin{equation}
   p(z)=\left( \begin{array}{cc}0_{\nu} & 0 \\ p_{21}(z) & p_{22}(z)\end{array}\right),
   \; \; \; \; \; \; \; \; q(z)=\left( \begin{array}{cc}
   q_{11}(z) & q_{12}(z) \\ q_{21}(z) & q_{22}(z)\end{array}\right)
\label{4.7}
\end{equation}
   and in view of (\ref{4.6}), $\; rank \ (q_{11}(z), \; q_{12}(z))=m$ at almost all 
   $z\in\C_+$. Multiplying $(q_{11}(z), \; q_{12}(z))$ by an appropriate unitary
   matrix $U$ on the right we obtain
   $$
   (q_{11}(z), \; q_{12}(z)) \ U=(\widetilde{q}_{11}(z), \; \widetilde{q}_{12}(z)),
   \; \; \; \; \; \; \det \; \widetilde{q}_{11}(z)\not\equiv 0.
   $$  
The pair $\{p, \; q\}$ is
   equivalent to  the pair $\{p_1, \; q_1\}$ defined as
   $$
   \left( \begin{array}{c}p_1(z) \\ q_1(z)\end{array}\right)=
   \left( \begin{array}{c}p(z) \\ q(z)\end{array}\right) \ U\Phi(z)\quad\mbox{where}\quad
\Phi(z)=\left( \begin{array}{cc}
   \widetilde{q}_{11}^{-1}(z) & -\widetilde{q}_{11}^{-1}(z)\widetilde{q}_{12}(z) \\
   0 & I_{m-\nu} \end{array}\right).
   $$
It follows from (\ref{4.7}) that the functions $p_1$ and $q_1$ are of the form
\begin{equation}
   p_1(z)=\left( \begin{array}{cc}0_{\nu} & 0 \\ \widetilde{p}_1(z) & \widetilde{p}(z)
   \end{array}\right) \; \; \; \; \; \; \; \; 
   q_1(z)=\left( \begin{array}{cc}I_{\nu} & 0 \\ \widetilde{q}_1(z) & \widetilde{q}(z)  
   \end{array}\right) 
\label{4.8}
\end{equation}
   and it remains to show that $\{p_1, \; q_1\}$ is equivalent to the pair defined
   in (\ref{4.65}). Indeed, (\ref{4.8}) implies that 
   $\{\widetilde{p}, \; \widetilde{q}\}\in {\bf N}_{m-\nu}$ and therefore, 
   $\det \; (\widetilde{p}(z)+i\widetilde{q}(z))\not\equiv 0$. 
   Substituting the pair (\ref{4.8}) into (\ref{4.5}) gives 
\begin{eqnarray*}
   S(z) & =&(p_1(z)-iq_1(z))(p_1(z)+iq_1(z))^{-1} \\ 
   &=&\left( \begin{array}{cc}-iI & 0 \\ \widetilde{p}_1-i\widetilde{q}_1 & 
\widetilde{p}-i\widetilde{q}
   \end{array}\right)\left( \begin{array}{cc}iI & 0 \\ \widetilde{p}_1+i\widetilde{q}_1 & 
   \widetilde{p}+i\widetilde{q}\end{array}\right)^{-1}\\ 
   &=&\left( \begin{array}{cc}-I & 0 \\
    i(\widetilde{p}-i\widetilde{q})(\widetilde{p}+i\widetilde{q})^{-1}(\widetilde{p}_1+i
\widetilde{q}_1) -   i(\widetilde{p}_1-i\widetilde{q}_1)& 
(\widetilde{p}-i\widetilde{q})(\widetilde{p}+i\widetilde{q})^{-1}
   \end{array}\right)\\ 
   & =&\left( \begin{array}{cc}-I & 0 \\ 0 & 
(\widetilde{p}-i\widetilde{q})(\widetilde{p}+i\widetilde{q})^{-1}
   \end{array}\right)
\end{eqnarray*}
   (to obtain the last equality we used the following: if the function $S={\scriptsize 
   \left( \begin{array}{cc}s_1 & 0 \\ s_2 & -I \end{array}\right)}$ is contractive 
   valued, then $s_2\equiv 0$). It is easily seen that the pair (\ref{4.65}) being substituted into 
(\ref{4.5}),
   leads to the same function $S$. By Lemma 4.2, the pairs (\ref{4.5}) and (\ref{4.8}) are 
equivalent.\qed     
\begin{La}
Let $R\in \C^{l\times 2m}$ be a 
   $J$--neutral matrix (i.e. $RJR=0$) and let $\; {\rm rank} \; R=\nu\leq 
   {\rm min}(m, \ l)$. Then there exist a $J$--unitary matrix $\Psi$ and an 
   invertible $T$ such that
\begin{equation}
   TR\Psi=\left( \begin{array}{cc}I_{\nu} & 0 \\ 0 & 0 \end{array}\right).
\label{4.9}
\end{equation}
\label{L:4.4}
\end{La}
{\bf Proof:} Since $\; {\rm rank} \; R=\nu$, there exists an
   invertible matrix $T$ such that
\begin{equation}
   TR=\left( \begin{array}{l}\widehat{R} \\ 0_{(m-\nu)\times 2m}\end{array}\right) 
\label{4.10}
\end{equation}
   where $\widehat{R}$ is a full rank $J$--neutral matrix. Let us endow the space
   $\C^{1\times 2m}$ with the indefinite inner product $[x, \ y]=yJx^*$. By
   (\ref{1.9}), the subspace   
   $$
   {\mathcal G}=\left\{g\in\C^{1\times 2m}: \; \; g=(\widehat{g}, \; 0) \; \; 
{\rm for \; \; 
   some} \; \; \widehat{g}\in\C^{1\times \nu}\right\}
   $$
   is $J$--neutral. The subspace ${\mathcal 
F}=\left\{f\in\C^{1\times 2m}: 
    \; f=\widehat{g}\widehat{R}, \; \widehat{g}\in\C^{1\times \nu}\right\}$ 
 $J$--neutral as well. Let us introduce the operator $\; {\bf \widehat{\Psi}} \; : \; 
   {\mathcal F}\rightarrow {\mathcal G} \; \;$ by
   $\; \; \widehat{g}\widehat{R} \ {\bf \widehat{\Psi}}=(\widehat{g}, \; 0)$. 
   Since ${\mathcal F}$ and ${\mathcal G}$ are $J$--neutral and 
   $\; \dim \; {\mathcal F}=\dim \; {\mathcal G}, \;$ the operator ${\bf 
   \widehat{\Psi}}$ is $J$--isometric and has equal defect numbers. Furthermore,
   ${\bf \widehat{\Psi}}$ is invertible and therefore,
   it admits a $J$--unitary extension ${\bf \Psi}$ to all of  $\C^{1\times 2m}$ 
   (\cite{ioh}). The matrix $\Psi$ of this extended operator in the standard basis is
   $J$--unitary  and satisfies $\; \widehat{R}\Psi=(I_{\nu}, \; 0) \;$ which both
   with (\ref{4.10}) implies (\ref{4.9}).\qed \vspace{3mm} \\
\begin{Rk}
{\rm Let $R=(R_1, \; R_2)\in \C^{l\times 2m}$ be a
   $J$--neutral matrix: $\; R_1R_2^*- R_2R_1^*=0$. Then 
   $\; {\rm rank} \; R={\rm rank} \; (R_1+iR_2)$. Indeed,
   $$
   {\rm rank} \; (R_1+iR_2)={\rm rank} \; (R_1+iR_2)(R_1+iR_2)^* =
   {\rm rank} \; (R_1R_1^*+R_2R_2^*)={\rm rank} \; RR^* ={\rm rank} \; R.
   $$}
\label{R:4.5}
\end{Rk}
   The following theorem is the degenerate analogue of Theorem \ref{T:1.3}.
\begin{Tm}
Let the Pick matrix $K_{n}$
   of the ${\bf HMP}$ be in the class $\widetilde{\mathcal H}_{m,n}$
   and let $\Theta$ be the $\C^{2m\times  2m}$--valued function defined by $(\ref{3.25})$.
   Then, there exists a J-unitary matrix $\Psi\in \C^{2m\times 2m}$
   such that  
\begin{enumerate}
\item All the functions $w\in {\mathcal R}(K_{n})$ are obtained 
by the formula
\begin{equation}
   w(z)=(a_{11}(z)p(z)+a_{12}(z)q(z))
   (a_{21}(z)p(z)+a_{22}(z)q(z))^{-1}
\label{4.11}  
\end{equation}   
   with the coefficient matrix $A(z)=(a_{ij}(z))=\Theta(z)\Psi\in {\bf W}$
   when the parameter $\{p,q\}$ varies in the set of all
   Nevanlinna pairs of the form
\begin{equation}
   \{p(z), \; q(z)\}=\left \{{\scriptsize \left( \begin{array}{cc}
   0_{\nu} & 0 \\ 0 & \widetilde{p}(z) \end{array} \right)} ,\hspace{2mm} 
   {\scriptsize \left( \begin{array}{cc} I_{\nu} &
   0 \\ 0 & \widetilde{q}(z) \end{array}  \right)} \right \}
\label{dop9}
\end{equation}
where $\{\widetilde{p}, \; \widetilde{q}\}\in {\bf N}_{m-\mu}$ and $\nu$ is the integer given 
by
$$
\nu=rank \; \left \{ (I_{m}, is_{0}, \ldots , is_{n-1})
   {\bf P}_{KerK_{n}}\right \}.
$$
\item Two pairs lead to the same function $w$ if and
   only if they are equivalent.
\end{enumerate}
\label{T:4.6}
\end{Tm}
   {\bf Proof:}  According to Theorem 1.1 the set ${\mathcal R}(K_{n})$
   coincides with the set of all solutions to the inequality (\ref{1.5})
   which is equivalent, by Lemma \ref{L:2.4}, to the following system
\begin{equation}
   \frac{w(z)-w(z)^{*}}{z-\bar{z}}-(Uw(z)+M)^{*}(I-zF)^{-*}K^{[-1]}
   (I-zF)^{-1}(Uw(z)+M)\geq 0, 
\label{4.12}  
\end{equation}   
\begin{equation}
   {\bf P}_{KerK}(I-zF)^{-1}\{Uw (z)+M\} \equiv 0.
\label{4.13}  
\end{equation}   
   It is easily seen that (\ref{4.12}) can be written as
   $$
   (w(z)^*, \; I) \left\{ \frac{J}{i(\bar{z}-z)}-\singcolb{U^*}{M^*}
   (I-zF^{*})^{-1} K^{[-1]} (I-zF)^{-1} (U, \; M) \right\}
   \singcolb{w(z)}{I} \geq 0
   $$
   and is equivalent, in view of (\ref{3.27}), to the inequality (\ref{4.1}) with the 
   function $\; \Theta \;$ defined by (\ref{3.25}) which is of the class
   ${\bf W}$ by Lemma \ref{L:3.3}. According to Theorem 4.1, all solutions 
   $\; w \;$ to the inequality (\ref{4.12}) are  parametrized by the linear
   fractional transformation (\ref{4.2}) when the parameter $\{p, \;  q\}$ varies 
   in the set ${\bf N}_{m}$ of all Nevanlinna pairs and satisfies (\ref{4.3}). 
   It remains to choose among these solutions all functions $\; w \;$ which
   satisfy also identity (\ref{4.13}). The rest of the proof is broken 
   into four steps which we now specify.

\medskip

\noindent
   {\bf Step 1:} The function $w(z)$ of the form (\ref{4.2}) stisfies the
   identity (\ref{4.13}) if and only if the corresponding parameter
   $\{p, \; q\}$ satisfies 
\begin{equation}
   {\bf P}_{KerK}\{Up(z)+Mq(z)\}\equiv 0.
\label{4.14}  
\end{equation}   
   {\bf Step 2:} If a pair $\{p, \; q\}\in {\bf N}_{m}$ satisfies (\ref{4.14})
   then it also satisfies (\ref{4.3}).\vspace{2mm} \\
   {\bf Step 3:} If a pair $\{p, \; q\}\in {\bf N}_{m}$ satisfies (\ref{4.14})
   then it is equivalent to some pair $\{p_{1}, \; q_{1}\}$ of the form
\begin{equation}
   \left( \begin{array}{c} p_{1}(z) \\ q_{1}(z)\end{array} \right)=\Psi
   {\scriptsize \left( \begin{array}{cc} 0_{\nu} & 0 \\ 0 & \widetilde{p}(z) \\
   I_{\nu} & 0 \\ 0 & \widetilde{q}(z)\end{array} \right)}  \sim
   \left( \begin{array}{c} p(z) \\ q(z) \end{array} \right)
\label{4.15}  
\end{equation}   
   for some $J$--unitary matrix $\; \; \Psi\in \C^{2m\times 2m} \; \;$ which 
   depends only on $\hspace{2mm}K_n\hspace{2mm}$ and a pair \\
   $\{\widetilde{p}, \; \widetilde{q}\}\in {\bf N}_{m-\nu}$, where
   $\; \nu=rank \; {\bf P}_{KerK}(U,M)=rank \; {\bf P}_{KerK}(U+iM).$

\medskip

   {\bf Proof of Step 1:} Let $\Theta=(\theta_{ij})$ be the function defined
   by (\ref{3.25}) and let $w$ be a function of the form (\ref{4.2}) for some
   pair $\{p, \; q\}\in {\bf N}_{m}$ which satisfies (\ref{4.3}). Then
   $$
   \left( \begin{array}{c} w(z) \\ I \end{array} \right) = \Theta(z)\left(
   \begin{array}{c} p(z) \\ q(z)\end{array} \right)  
   (\theta_{21}(z)p(z)+\theta_{22}(z)q(z))^{-1}
   $$
   and therefore, identity (\ref{4.13}) is equivalent to
\begin{equation}
   {\bf P}_{KerK}(I-zF)^{-1}(U,\hspace{1mm}M)\Theta(z) \left( \begin{array}{c} p(z) \\
   q(z) \end{array} \right) \equiv 0.
\label{4.16}  
\end{equation}   
   Using (\ref{1.8}), (\ref{3.25}) and the identity 
   $$
   K(I-zF^{*})^{-1}-(I-zF)^{-1}K=z(I-zF)^{-1}(KF^{*}-FK)(I-zF^{*})^{-1}
   $$
   we get
$$
   (I-zF)^{-1}(U, \; M)\Theta(z)=K(I-zF^{*})^{-1}K^{[-1]}(U, \; M)
   +(I-zF)^{-1}\left\{I-KK^{[-1]}\right\}(U, \; M). 
$$
Substituting the latter equality into (\ref{4.16}) gives 
$$
   {\bf P}_{KerK}(I-zF)^{-1}\{I-KK^{[-1]}\}\left(Up(z)+Mq(z)\right)\equiv 0
$$
which on account of (\ref{2.4}), can be written as
\begin{equation}
   \{I+z{\bf P}_{KerK}F(I-zF)^{-1}(I-KK^{[-1]})\}{\bf P}_{KerK}\left(Up(z)+Mq(z)\right)\equiv 0.
\label{4.18}  
\end{equation}   
   Since the matrix $\{I+z{\bf P}_{kerK}F(I-zF)^{-1}(I-KK^{[-1]})\}$  
   is nondegenerate, (\ref{4.18}) implies (\ref{4.14}). 

\medskip

   {\bf Proof of Step 2:} Let a pair $\{p,\hspace{1mm}q\}\in {\bf N}_{m}$ satisfy the
   condition (\ref{4.14}). We introduce the pair
\begin{equation}
   \left( \begin{array}{c} p_{0}(z) \\ q_{0}(z)\end{array} \right) =
   \Theta(z) \left( \begin{array}{c} p(z) \\ q(z)\end{array} \right)
\label{4.19}  
\end{equation}   
   and show that $\; det \ q_{0}(z)\not\equiv 0$. Indeed, suppose
   that the point $\lambda\in \C_{+}$ and the nonzero vector ${\bf h}\in \C^{m}$ 
   are such that $\; \det \ \Theta(\lambda)\neq 0 \;$ and 
\begin{equation}
   q_{0}(\lambda){\bf h}=0.
\label{4.20}  
\end{equation}   
   Since ${\bf h}^{*}\left(p(\lambda)^{*},\hspace{1mm}q(\lambda)^{*}\right)
   \Theta(\lambda)^{*}J\Theta(\lambda)
   \left( \begin{array}{c} p(\lambda) \\ q(\lambda)\end{array} \right){\bf h}
   ={\bf h}^{*}\left (p_{0}(\lambda)^{*},\hspace{1mm}0\right )J \left( \begin{array}{c}
   p_{0}(\lambda) \\ 0\end{array} \right){\bf h}=0$, then
   $$
   0  \leq {\bf h}^{*} \left (p(\lambda)^{*},\hspace{1mm}q(\lambda)^{*}\right) J \left(
   \begin{array}{c} p(\lambda) \\ q(\lambda)\end{array} \right) {\bf h}  
     ={\bf h}^{*} \left (p(\lambda)^{*},\hspace{1mm}q(\lambda)^{*}\right )
   \left \{J-\Theta(\lambda)^{*}J\Theta(\lambda)\right \}
   \left( \begin{array}{c} p(\lambda) \\ q(\lambda)\end{array}\right){\bf h},
   $$
due to (\ref{1.10}). Substituting (\ref{3.26}) into this last inequality leads us to
\begin{equation}
   K(I-\lambda
   F^{*})^{-1}K^{[-1]}\{Up(\lambda)+Mq(\lambda)\}{\bf h}=0.
\label{4.21}     
\end{equation}      
   It follows from (\ref{3.25}) and (\ref{4.19}) that
\begin{equation}
   p_{0}(\lambda)=p(\lambda)+\lambda M^{*}(I-\lambda F^{*})^{-1}K^{[-1]}
   \{Up(\lambda)+Mq(\lambda)\}.
\label{4.22}     
\end{equation}      
   Since $M=FKU$ (see (\ref{1.7})), then $\hspace{2mm}\lambda M^{*}(I-\lambda F^{*})^{-1}
   =U^{*}K(I-\lambda F^{*})^{-1}-U^{*}K.\hspace{2mm}$
   Substituting this last equality into (\ref{4.22}) and taking into account
   (\ref{2.4}), (\ref{4.14}), (\ref{4.21}) and the evident equalities $U^{*}U=I_{m}$
   and $U^{*}M=0$ we receive 
   $$
   p_{0}(\lambda){\bf h}=p(\lambda){\bf h}-U^{*}KK^{[-1]}\{Up(\lambda)+Mq(\lambda)\}
   {\bf h} +U^{*}K(I-zF^{*})^{-1}K^{[-1]} \{Up(\lambda)+Mq(\lambda)\}{\bf h}
   $$
   $$\hspace{-16mm}
   =U^{*}(I-KK^{[-1]})\{Up(\lambda)+Mq(\lambda)\}{\bf h}+(I-UU^{*})p(\lambda)
   -U^{*}Mq(\lambda) 
   $$
   $$ \hspace{-49mm}
   =U^{*}(I-KK^{[-1]}){\bf P}_{KerK}\{Up(\lambda)+Mq(\lambda)\}{\bf h}=0.
   $$
   Since $\hspace{2mm}det \; \Theta(\lambda)\neq 0, \;$ the equality
   $\; \; p_{0}(\lambda)h=0 \; \;$ both with (\ref{4.19}) and (\ref{4.20}) implies
   $$ 
   \left( \begin{array}{c} p(\lambda) \\ q(\lambda)\end{array} \right){\bf h}=
   \Theta(\lambda)^{-1} \left( \begin{array}{c} p_{0}(\lambda) \\
   q_{0}(\lambda)\end{array} \right){\bf h}=0
   $$
and since $\lambda$ is a arbitrary point, the latter equality contradicts to 
   the nondegeneracy of the pair $\{p, \; q\}$. 

\medskip

   {\bf Proof of Step 3:} Using (\ref{1.8}) we obtain that the matrix
   $\; {\bf P}_{KerK}(U, \; M) \;$ is $J$--neutral:
   $$
   {\bf P}_{KerK}(U, \; M)J \left( \begin{array}{c} U^{*} \\ M^{*}\end{array}
   \right) {\bf P}_{KerK}=i{\bf P}_{KerK}(KF^{*}-FK){\bf P}_{KerK}=0.
   $$
Thus, by Remark 4.5,
   $$
   \mu={\rm rank} \; ({\bf P}_{KerK}(U, \; M))={\rm rank} \; ({\bf P}_{KerK}(U+iM))=
   {\rm rank} \; \left\{(I_{m},is_{0},\ldots,is_{n-1}){\bf P}_{KerK_{n}}\right\}.
   $$
   According to Lemma 4.4, there exist a $J$--unitary matrix $\Psi$ and
an invertible $T$ such that 
\begin{equation}
   T{\bf P}_{KerK}(U, \; M)\Psi=\left( \begin{array}{cc}I_{\nu} & 0 \\ 0 & 0 
   \end{array}\right).
\label{4.23}     
\end{equation}      
   Let $\{p_{2}, \; q_{2}\}$ be the pair defined by
\begin{equation}
   \left( \begin{array}{c}p(z) \\ q(z)\end{array}\right)=
   \Psi \; \left( \begin{array}{c}p_2(z) \\ q_2(z)\end{array}\right).
\label{4.24}     
\end{equation}      
   On account of (\ref{4.23}) and (\ref{4.24}), condition (\ref{4.3}) can be rewritten as 
   $\; \left(I_{\nu}, \; 0\right) \ p_2(z)\equiv 0 \;$
   and by  Lemma 4.3, the pair
   $\{p_2, \; q_2\}$ is equivalent to some pair of the form (\ref{4.65}), i.e.,
   $$ 
   \left( \begin{array}{c} p(z) \\ q(z) \end{array} \right)=
   \Psi \; \left( \begin{array}{c}p_2(z) \\ q_2(z)\end{array}\right)
   \sim \Psi {\scriptsize \left( \begin{array}{cc} 0_{\nu} & 0 \\ 0 & \widehat{p}(z) \\ 
   I_{\nu} & 0 \\ 0 & \widehat{q}(z)\end{array} \right)}=
   \left( \begin{array}{c} p_{1}(z) \\ q_{1}(z)\end{array} \right)   
   $$
   which completes the proof of Step 3.

\medskip

   Substituting (\ref{4.15}) into (\ref{4.2}) and taking into account that the
   equivalent pairs lead under the linear fractional transformation to
   the same function $w(z)$, we finish the proof  of theorem.\qed

\medskip

By Remark \ref{R:2.9}, the condition $K_{n}\in\widetilde{\mathcal H}_{m,n}$ 
   is not restrictive and hence, the received in Theorem 3.3 description is 
   applicable to the general situation $K_{n} \in {\bf {\mathcal 
H}}_{m,n}$.

   \section{Correction of erratae in \cite{bol}}
\setcounter{equation}{0}

The following result was formulated in \cite{bol} (see Lemmas 2.5, 2.10 and 2.11 there).
\begin{La}
 Let  $K_{n}=(s_{i+j})_{i,j=0}^n \in {\bf {\mathcal H}}_{m,n}$ and let ${\mathcal L}$ be the
subspace of $\C^{1 \times m}$ given in \eqref{2.6}. The following are equivalent:
\begin{enumerate}
\item $K_n$ admits a positive semidefinite Hankel extension.
\item ${\bf P}_{Ker K_{n-1}} {\scriptsize \left( \begin{array}{c}
   s_{n+1} \\ \vdots \\ s_{2n} \end{array} \right)} =0$.
\item The block  $s_{2n}$ is of the form
\begin{equation}
   s_{2n}=(s_{n},\ldots,s_{2n-1})K_{n-1}^{[-1]}(s_{n},\ldots,s_{2n-1})^{*}+R
\label{5.1}
\end{equation}
   for some positive semidefinite matrix $R\in \C^{m \times m}$
   which vanishes on the subspace ${\mathcal L}$ and does not depend on
the choice of $K_{n-1}^{[-1]}$.  
\item The associated truncated Hamburger moment problem admits an ``exact'' solution $\sigma$ 
such
that
$$
   \int^{\infty}_{-\infty}\lambda^{k}d\sigma(\lambda)=s_{k} \hspace{10mm}
   (k=0,\ldots,2n).
$$
\end{enumerate}
\label{L:5.1}
\end{La}
The proofs of implications $(1)\Rightarrow (4)\Rightarrow (2)\Leftrightarrow (3)$ presented in 
\cite{bol} are correct; they are reproduced in Lemmas \ref{2.5} and \ref{2.6} above. To complete 
the proof, it suffices to justify $(2)\Rightarrow (1)$, that is, in our current terminology, to 
show that 
\begin{equation}
\widetilde{\mathcal H}_{m,n}\subseteq {\mathcal H}_{m,n}^+.
\label{2.13b} 
\end{equation}
This inclusion together with \eqref{2.13a} implies that all three classes introduced in Section 
2 coincide.

\bigskip

{\bf \underline{Proof of \eqref{2.13b}}:} Let $K_n\in\widetilde{\mathcal H}_{m,n}$. Plug in the 
Nevanlinna pair $\{p,q\}=\{0_m,  I_m\}$ (which is certainly of the form \eqref{dop9})
into formula \eqref{4.11} to get a solution $w(z)=a_{12}(z)a_{22}(z)^{-1}$ from ${\mathcal 
R}(K_n)$. This Pick function $w$ is rational (since $A$ is) and takes Hermitian values 
at every real point at which it is analytic (since $A$ is $J$-unitary on $\R$). Then the measure 
$\sigma$ from the Herglotz representation \eqref{1.4} of $w$ is finitely atomic and therefore,
the integrals $\int_{-\infty}^\infty \lambda^N d\sigma(\lambda)$ exists for every $N\ge 0$.
Since this measure solves the associated ${\bf HMP}$, it satisfies \eqref{1.2} and \ref{1.3}.
By virtue of \eqref{2.19}, the Hankel block matrix 
$$
\widetilde{K}_n=\left( \begin{array}{cc}K_{n-1} & \begin{array}{cc} s_n & s_{n+1} \\
\vdots & \vdots \\ s_{2n-1} & s_{2n} \end{array} \\
 \begin{array}{ccc} s_n & \cdots & s_{2n-1} \\ s_{n+1} & \cdots & s_{2n}\end{array} &
\begin{array}{cc}s & s_{2n+1} \\ s_{2n+1} & s_{2n+2}\end{array}\end{array}\right)
$$
is positive semidefinite, where we have set 
$$
s=\int_{-\infty}^\infty\lambda^{2n}d\sigma(\lambda), \quad 
s_{2n+1}=\int_{-\infty}^\infty\lambda^{2n+1}d\sigma(\lambda),\quad
s_{2n+2}=\int_{-\infty}^\infty\lambda^{2n+2}d\sigma(\lambda).
$$
The Hankel block matrix $K_{n+1}:=(s_{i+j})_{i,j=0}^{n+1}$ extends $K_n$ and is positive 
semidefinie. Indeed, by \eqref{1.3}, we have $K_{n+1}\ge \widetilde{K}_n\ge 0$. Thus 
$K_n\in{\mathcal H}_{m,n}^+$ which completes the proof.\qed
\begin{Rk}
{\rm The proof of implication $(2)\Rightarrow (1)$ presented in \cite{bol} does not rely on 
interpolation Theorem \ref{T:4.6}. The extending matrices $s_{2n+1}$ and $s_{2n+2}$ were
constructed directly in terms of the given $s_0, \ldots, s_{2n}$. Unfortunately, the 
construction turned out to be wrong. The author was very glad to learn that 
correct explicit proofs of the above implication have been recently obtained \cite{dyuk, wo}.}
\label{R:5.2}
\end{Rk}

\end{document}